\documentclass{amsart}

\usepackage{amsmath,amssymb,amsthm}
\usepackage{graphicx,tikz}
\newtheorem{theorem}{Theorem}
\newtheorem*{thm}{Theorem}

\theoremstyle{definition}

\theoremstyle{remark}

\begin{document}

\title[]{On the Regularization Effect of Stochastic Gradient Descent applied to Least Squares}
\subjclass[2010]{65F10, 65K10, 65K15, 90C06, 93E24} 
\keywords{Stochastic Gradient Descent, Kaczmarz method, Least Squares}
\thanks{S.S. is supported by the NSF (DMS-1763179) and the Alfred P. Sloan Foundation.}

\author[]{Stefan Steinerberger}
\address{Department of Mathematics, University of Washington, Seattle}
\email{steinerb@uw.edu}

\begin{abstract} We study the behavior of stochastic gradient descent applied to $\|Ax -b \|_2^2 \rightarrow \min$ for invertible $A \in \mathbb{R}^{n \times n}$. We show that there is an explicit constant $c_{A}$ depending (mildly) on $A$ such that
$$ \mathbb{E} ~\left\| Ax_{k+1}-b\right\|^2_{2} \leq \left(1 + \frac{c_{A}}{\|A\|_F^2}\right) \left\|A x_k -b \right\|^2_{2} - \frac{2}{\|A\|_F^2} \left\|A^T A (x_k - x)\right\|^2_{2}.$$
This is a curious inequality:  the last term has one more matrix applied to the residual $u_k - u$ than the remaining terms: if $x_k - x$ is mainly comprised of large singular vectors, stochastic gradient descent leads to a quick regularization. For symmetric matrices, this inequality has an extension to higher-order Sobolev spaces. This explains a (known) regularization phenomenon: an energy cascade from large singular values to small singular values smoothes. 
 \end{abstract}

\maketitle

\section{Introduction}
\subsection{Stochastic Gradient Descent.} In this paper, we consider the finite-dimensional linear inverse problem
$$ Ax =b,$$
where $A \in \mathbb{R}^{n \times n}$ is an invertible matrix, $x \in \mathbb{R}^n$ is the (unknown) signal of interest and $b$ is a given right-hand side. Throughout this paper, we will use $a_1, \dots, a_n$ to denote the rows of $A$. Equivalently, we will try to solve the problem 
$$\|Ax - b\|^2  = \sum_{i=1}^{n} \left( \left\langle a_i, x \right\rangle - b_i\right)^2 \rightarrow \min.$$
Following Needell, Srebro \& Ward \cite{need3}, we can interpret this problem as
$$ \sum_{i=1}^{n} f_i(x)^2 \rightarrow \min \qquad \mbox{where} \qquad f_i(x) = \left\langle a_i, x \right\rangle - b_i.$$
The Lipschitz constant of $f_i$ is $\|a_i\|_{\ell^2}$ which motivates the following basic form of stochastic gradient descent: pick one of the $n$ functions with likelihood proportional to the Lipschitz constant and then do a gradient descent on this much simpler function resulting in the update rule
$$ x_{k+1} = x_k - \frac{b_i - \left\langle a_i, x_k\right\rangle}{\|a_i\|^2}a_i.$$
This is also known as the \textit{Algebraic Reconstruction Technique} (ART) in computer tomography \cite{gordon, her, her2, natterer}, the \textit{Projection onto Convex Sets Method} \cite{cenker, deutsch, deutsch2, gal, sezan} and the \textit{Randomized Kaczmarz method} \cite{bai, eldar, elf, gordon2, gower, gower2, lee, leventhal, liu, ma, moor, need, need2, need25, need3, need4, nutini, popa, stein, stein2, strohmer, tan, zhang, zouz}.  Strohmer \& Vershynin \cite{strohmer} showed that
$$ \mathbb{E} \|x_k - x\|^2 \leq \left(1 - \frac{1}{\|A^{-1}\|^2 \|A\|_F^2}\right)^k \|x_0 - x\|^2,$$
where $\|A\|_F^2$ is the Frobenius norm. In practice, the algorithm often converges a lot faster initially and this was studied in \cite{jiao, jin, stein}. In particular, \cite{stein} obtains an identity in terms of the behavior with regards to the singular values showing that singular vectors associated to large singular values are expected to undergo a more rapid decay. Motivated by this insight, we provide rigorous bounds that quantify this energy cascade from large singular values to small singular values by identifying an interesting inequality for SGD when applied to Least Squares.

\subsection{A Motivating Example.} We discuss a simple example that exemplifies the phenomenon that we are interested in. Let us take $A \in \mathbb{R}^{100 \times 100}$ by picking each entry independently at random from $\mathcal{N}(0,1)$ and then normalizing the rows to $\|a_i\|=1$. The right-hand side is $b=(1,1,\dots,1)$ and we initialize with $x_0 = \textbf{0} \in \mathbb{R}^n$.

\begin{center}
\begin{figure}[h!]
\begin{tikzpicture}[scale=1]
\node at (0,0) {\includegraphics[width=0.45\textwidth]{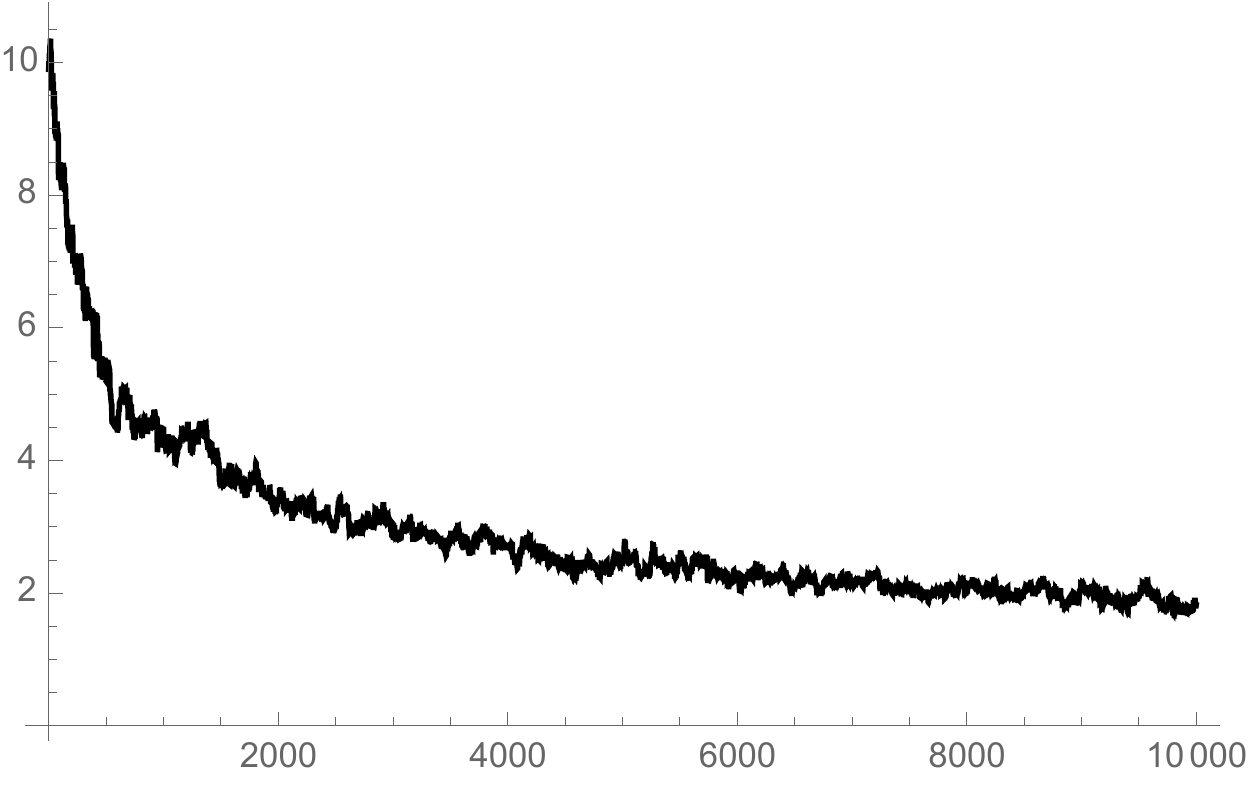}};
\node at (1,1) {$\|A x_k -b\|_{\ell^2}$};
\node at (6,0) {\includegraphics[width=0.45\textwidth]{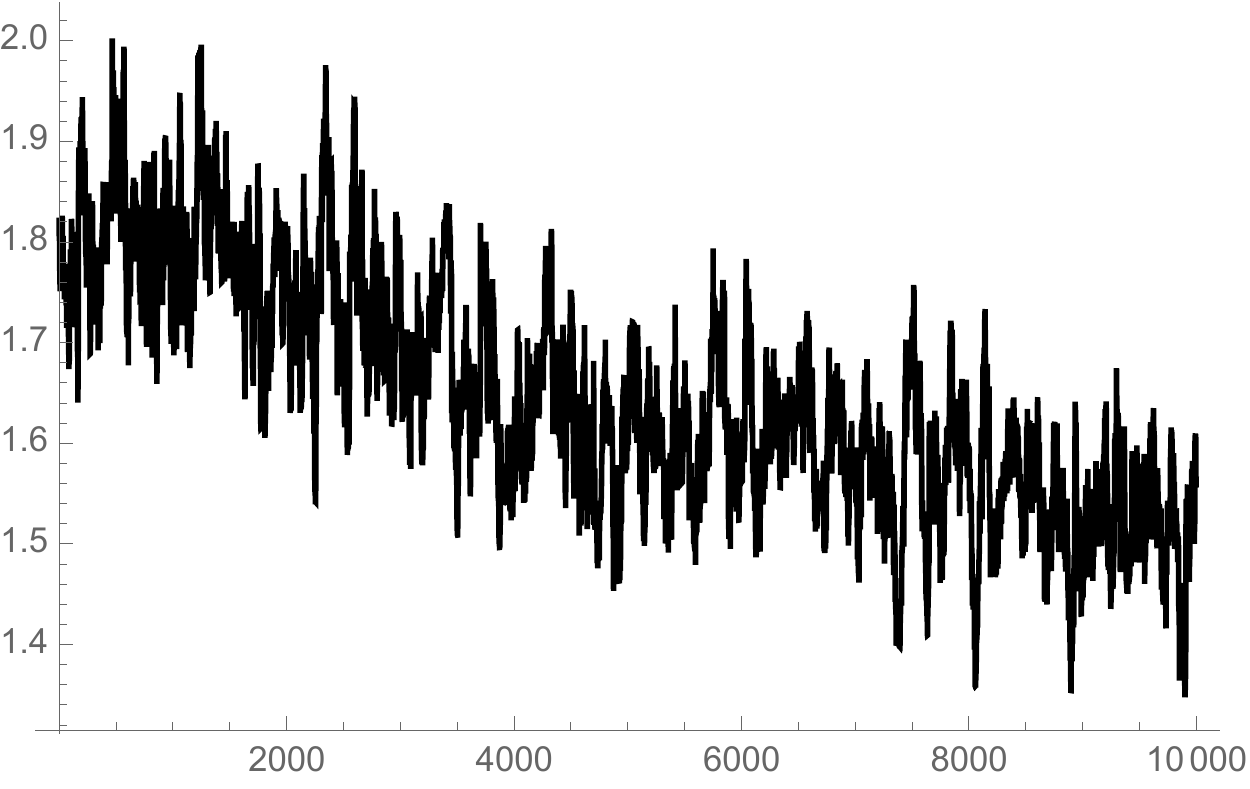}};
\node at (7,1) {$\|A x_k -b\|_{\ell^2}$};
\end{tikzpicture}
\caption{The size of $\|A x_k - b\|_{\ell^2}$ for $k=1, \dots, 10000$ (left) followed by the next $10000$ iterations $\|A x_{10000+k} - b\|_{\ell^2}$. We observe rapid initial decay which then slows down. }
\end{figure}
\end{center}

The picture tells a very interesting story: the error in $\|Ax_k - b\|$ decays initially quite rapidly before then stabilizing in a certain regime. Moreover, for the example shown in Figure 1, $\|x_0\|=0$ and $\| x_{20000}\| \sim 28$ which is not even close to the true solution $\|x\| \sim 128$ -- nonetheless, the approximation of $A x_k$ to $b$ is quite good.
This leaves us with a curious conundrum: we have a good approximation $x_k$ of the true solution in the sense that $A x_k \sim b$ even though $x_k$ is not very close to $x$. One way this can be achieved is if $x_k - x$ is mainly a linear combination of small singular vectors of $A$. This is related to the following result recently obtained by the author.
\begin{thm}[\cite{stein}] Let $v_{\ell}$ be a (right) singular vector of $A$ associated to the singular value $\sigma_{\ell}$. Then, for the sequence $(x_k)_{k-0}^{\infty}$ obtained in this randomized manner
$$\mathbb{E} \left\langle x_{k} - x, v_{\ell} \right\rangle = \left(1 - \frac{\sigma_{\ell}^2 }{\|A\|_F^2}\right)^k\left\langle x_0 - x, v_{\ell}\right\rangle.$$
\end{thm}

Here, $\|A\|_F$ denotes the Frobenius norm. This shows that we expect $x_k - x$ to be indeed mainly a linear combination of singular vectors associated to small singular values since those are the ones undergoing the slowest decay. It also mirrors the bound obtained by Strohmer \& Vershynin \cite{strohmer} since 
$$\sigma_{\ell} \geq \sigma_n = \|A^{-1}\|^{-1}.$$
 While being interesting in itself, this identity by itself does not fully explain the behavior shown above: it is only in expectation with no control of the variance. Moreover, the inner product does initially undergo some fluctuations. Taking the same type of matrix as above, we see an example of such fluctuations in Fig. 2.

\begin{center}
\begin{figure}[h!]
\begin{tikzpicture}[scale=1]
\node at (0,0) {\includegraphics[width=0.7\textwidth]{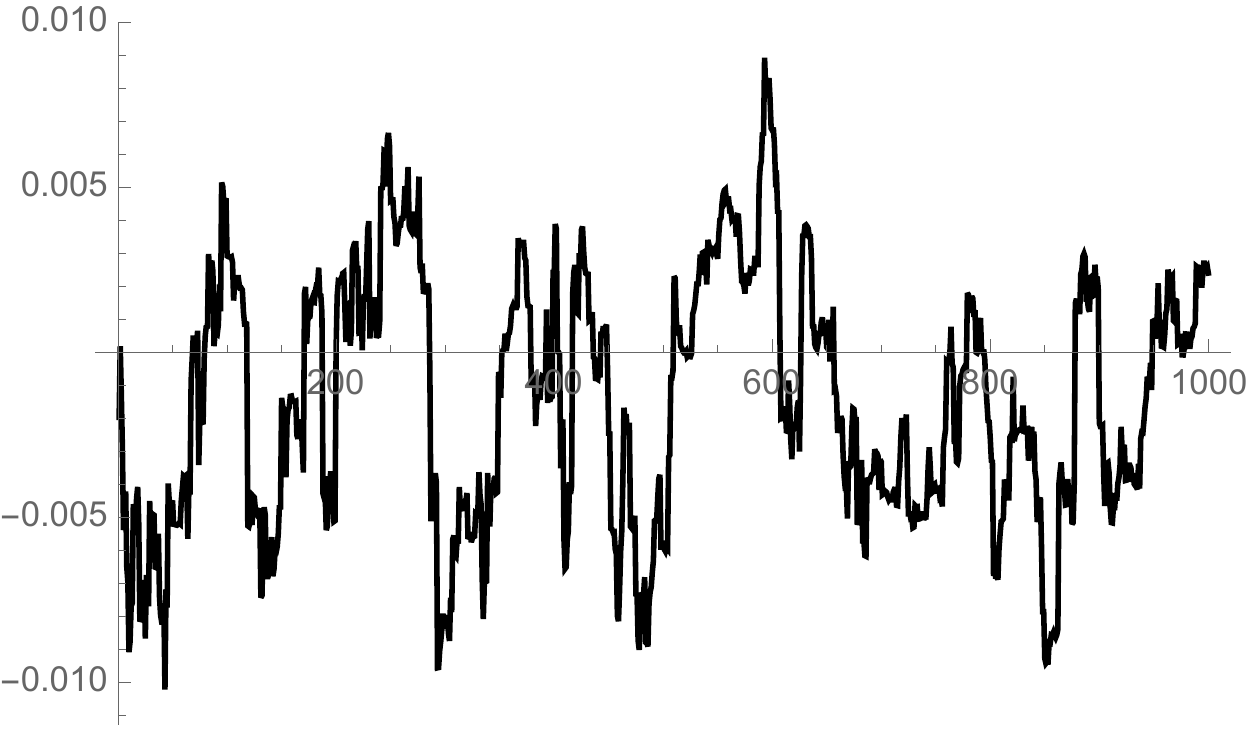}};
\node at (2.5,1.5) {\Large $\left\langle \frac{x_k -x}{\|x_k - x\|}, v_1\right\rangle$};
\end{tikzpicture}
\caption{The evolution of the normalized residual against the leading singular vector $v_1$: fluctuations around the mean.}
\end{figure}
\end{center}

\subsection{Related results.} This type of question is well studied, we refer to Ali, Dobridan \& Tibshirani \cite{ali}, Defossez \& Bach \cite{def}, 
Jain, Kakade, Kidambi, Netrapalli, Pillutla \& Sidford \cite{jain},
Neu \& Rosasco \cite{neu}, Oymak \& Soltanolkotabi \cite{oymak}, Schmidt, Le Roux \& Bach \cite{schm} and references therein. The connection of Stochastic Gradient Descent applied to Least Square Problems and the Randomized Kaczmarz Method has been pointed out by Needell, Srebro \& Ward \cite{need3}.
We also mention the papers by Jiao, Jin \& Lu \cite{jiao} and Jin \& Lu \cite{jin} who studied a similar question and noted that there is energy transfer from large singular values to small singular values.

\section{Results}
\subsection{Main Result.}
The main goal of this note is to provide a simple explanation for the rapid initial regularization: the expected decay of $\|A(x_k - x)\|$ under SGD can be bounded from above by a term involving $\left\| A^T A (x_k -  x)\right\|_2$: this is the same term except that the matrix has been applied to the existing quantity one more time. This increases the norm of the underlying vector except when $A(x_k - x)$ is mainly the linear combination of singular vectors with small singular values. So while this is not the case, we actually inherit strong decay properties and this leads to the rapid initial regularization. We now make this precise.

\begin{theorem}
Let $A \in \mathbb{R}^{n \times n}$ be invertible and consider $\|Ax - b\|^2_2 \rightarrow \min$ via the stochastic gradient descent method introduced above. Abbreviating
$$ \alpha = \max_{1 \leq i \leq n}{ \frac{\|A a_i\|^2}{\|a_i\|^2}},$$
we have
$$ \mathbb{E}~\|A x_{k+1}- b\|_2^2 \leq  \left(1 + \frac{\alpha}{\|A\|_F^2}\right)\|A x_k - b\|_2^2 - \frac{2}{\|A\|_F^2} \left\| A^T (A x_k -  b)\right\|_2^2.$$
\end{theorem}
The inequality also holds for $A \in \mathbb{R}^{m \times n}$ with $m \geq n$ as long as $Ax = b$ has a unique solution.
We note that $\alpha \geq 1$ is usually quite small: for random matrices, we would expect $\alpha \sim 2$. This is also true for matrices discretizing PDEs via finite elements.
The main point of the inequality is that the last term has an additional factor of $A^T$: we can rewrite it as
$$\left\| A^T (A x_k - b)\right\|_2^2 = \| A^T A (x_k - x)\|_2^2.$$
This shows that the presence of large singular vectors in $x_k - x$ forces large decay on $\|A(x_k - x)\|_2^2$. Conversely, once the algorithm
has reached the plateau phase (see Fig. 1), the terms
$$ \alpha \|A(x_k - x)\|_2^2 \qquad \mbox{and} \qquad \| A^T A (x_k - x)\|_2^2$$
are nearly comparable. Thus, this forces
 $x_k - x$  to be mostly orthogonal to most singular vectors corresponding to large singular values: that,
however, shows that it is mainly comprised of small singular vectors and thus explains why $\|A(x_k - x)\| \ll \|x_k - x\|$ is possible in cases where 
$x_k$ is far away from $x$. In particular, this suggests why the method could be effective for the problem of finding a vector $x$ so that $A x \approx b$.  One way is to initialize stochastic gradient descent with $x_0 = 0$ and run it for a while -- due to the difference in scales, second-order
norms regularizing first-order norms, we observe that $A x_k$ converges quite rapidly; whether it converges to something sufficiently close to $b$ for the purpose at hand, is a different question.

\subsection{A Sobolev Space Explanation.} An interesting way to illustrate the result is in terms of partial differential equation. Suppose we try to solve $-\Delta u = f$ on some domain $\Omega \subset \mathbb{R}^n$. After a suitably discretization, this results in a discrete linear system $Lu = f$, where $L \in \mathbb{R}^{n \times n}$ is a discretization of the Laplacian $-\Delta$. By an abuse of notation, $u$ denotes a decent approximation of the continuous solution and $f$ a discretization of the continuous right-hand side. However, we also have more information: since $L$ discretizes the Laplacian, we expect that
$$ \left\langle Lu, u\right\rangle \sim \int_{\Omega}{ |\nabla u|^2 dx} \qquad \mbox{and} \qquad \left\langle Lu, Lu \right\rangle \sim \int_{\Omega}{ |\Delta u|^2 dx}.$$
Here, the first term correspond to the size of $u$ in the Sobolev space $\dot H^1$ while the second term is the size of $u$ in the Sobolev space $\dot H^2$. In fact, this is a common way to define discretized approximations of Sobolev space, also known as the spectral definition since they are defined in terms of the spectrum of $L$. Suppose now we compute a sequence of approximations $u_k$ via the method outlined above. Then Theorem 1 can be rephrased as
$$ \mathbb{E} ~\| u_{k+1} - u\|^2_{\dot H^1} \leq   \left(1 + \frac{\alpha}{\|A\|_F^2}\right)\|u_k -u \|_{\dot H^1}^2 - \frac{2}{\|A\|_F^2} \left\| u_k - u\right\|_{\dot H^2}^2$$
What is of great interest here is that the decay of the error in $\dot H^1$ is driven by decay of the error in $\dot H^2$ (which is usually larger). 

\subsection{An Example.} We illustrate this with an example. Choosing $A \in \mathbb{R}^{500 \times 500}$ at random (and then, for convenience, normalize the rows to $\|a_i\|=1$), we solve $Ax = (1,1,\dots, 1)$ starting with a random initial vector $x_0$ where each entry is chosen independently from a standardized $\mathcal{N}(0,1)$ distribution. We consider both the evolution of $\|A x_k - b\|_{2}^2$ across multiple runs as well as the size of
$$  \frac{\alpha}{\|A\|_F^2}\|A x_k - b\|_2^2 - \frac{2}{\|A\|_F^2} \left\| A^T (A x_k -  b)\right\|_2^2,$$
which is the term from our Theorem quantifying the expected decay at each step.

\begin{center}
\begin{figure}[h!]
\begin{tikzpicture}[scale=1]
\node at (0,0) {\includegraphics[width=0.45\textwidth]{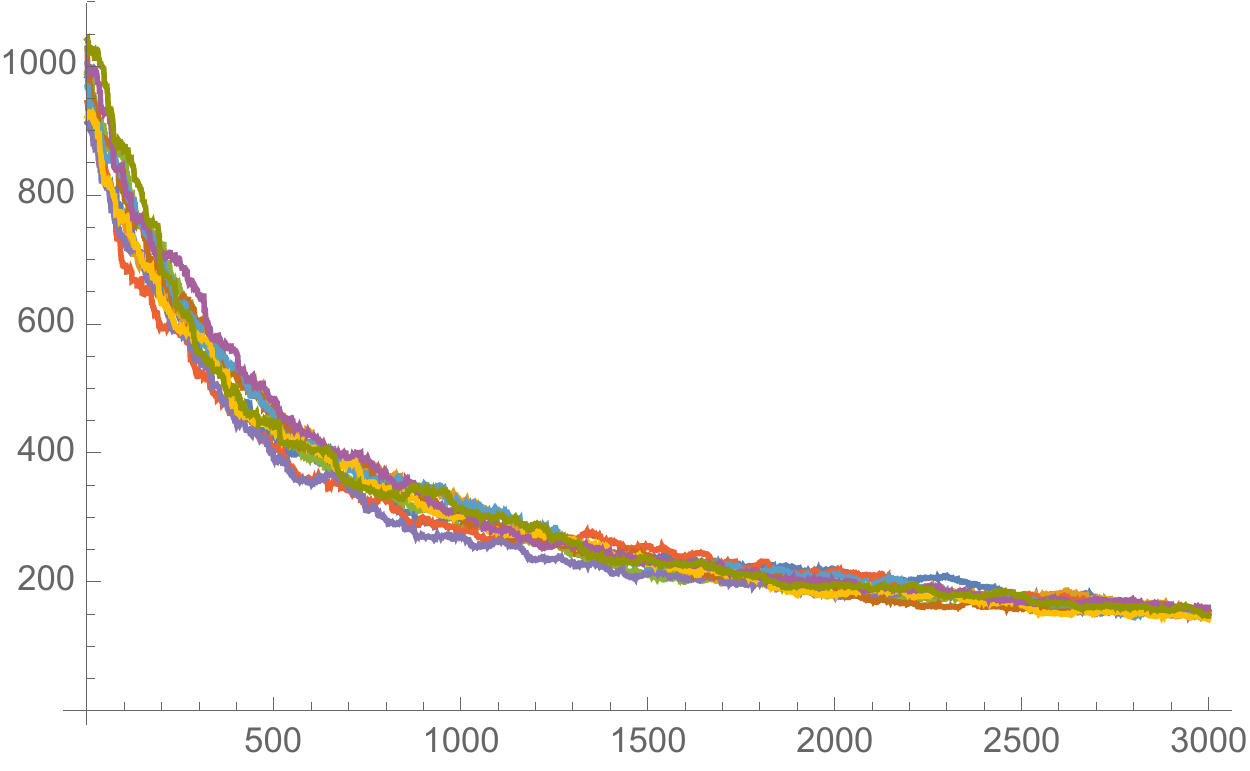}};
\node at (1,1) {$\|A x_k -b\|_{\ell^2}$};
\node at (6,0) {\includegraphics[width=0.45\textwidth]{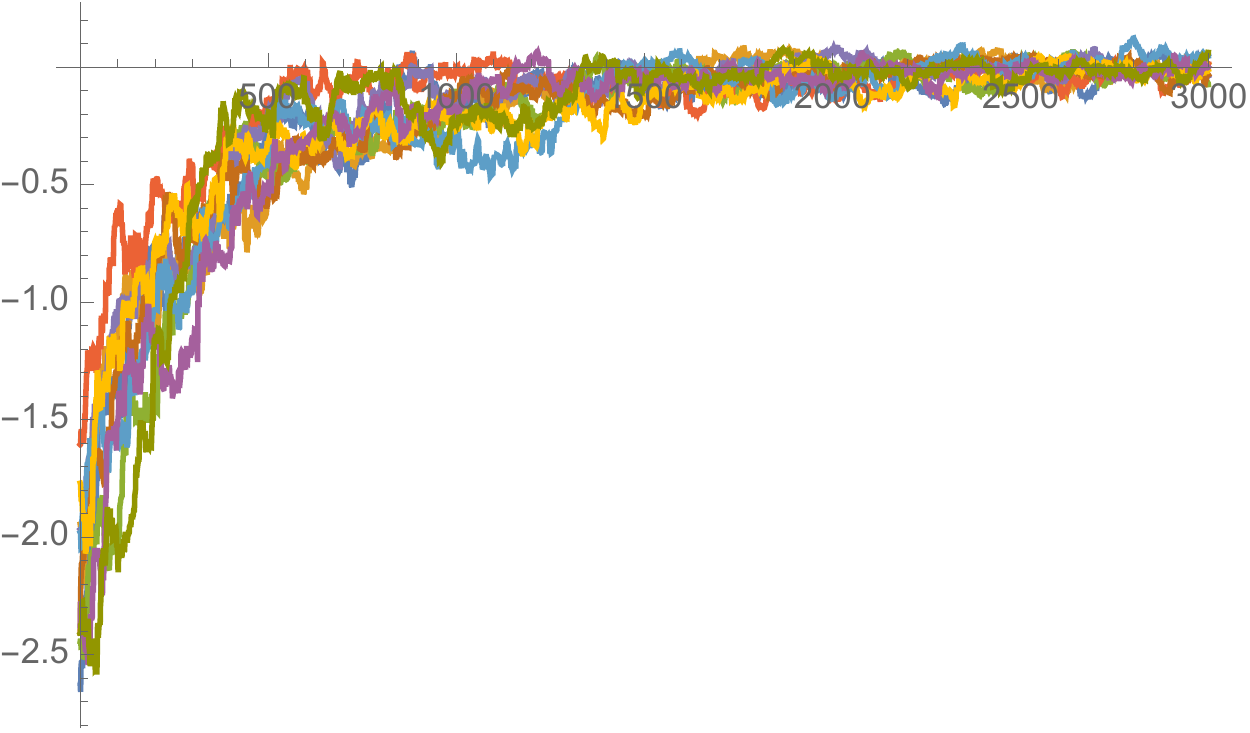}};
\node at (6.7,-0.2) {\tiny $  \frac{\alpha}{\|A\|_F^2}\|A x_k - b\|_2^2 - \frac{2}{\|A\|_F^2} \left\| A^T (A x_k -  b)\right\|_2^2
$};
\end{tikzpicture}
\caption{$\|A x_k - b\|_{\ell^2}$ for $k=1, \dots, 3000$ (left) and the decay guaranteed by Theorem in expectation (right) over multiple runs. }
\end{figure}
\end{center}

We see in Fig. 3. that over 3000 periods, the approximation decays roughly by a factor $\sim 800$ (with little variation across multiple runs). The bound in the Theorem implies an expected decay of $-0.23$ per time-step which, across 3000 time steps, leads to a total of roughly $\sim 696$ in decay.

\subsection{Higher Powers.} If the matrix $A \in \mathbb{R}^{n \times n}$ is symmetric, we can extend the result to higher powers of the matrix.
\begin{theorem}
Let $A \in \mathbb{R}^{n \times n}$ be symmetric and invertible. When solving $\|Ax - b\|^2_2 \rightarrow \min$ via the stochastic gradient descent method outlined above, we have, for any $\ell \in \mathbb{N}$ with
$$ \alpha_{\ell} = \max_{1 \leq i \leq n}{ \frac{\|A^{\ell} a_i\|^2}{\|a_i\|^2}},$$
the estimate
$$ \mathbb{E}~\|A^{\ell} (x_{k+1}- x^*)\|_2^2 \leq  \left(1 + \frac{\alpha_{\ell}}{\|A\|_F^2}\right)\|A^{\ell}( x_k - x^*)\|_2^2 - \frac{2}{\|A\|_F^2} \left\| A^{\ell+1} ( x_k -  x^*)\right\|_2^2.$$
\end{theorem}
This shows that the same phenomenon does indeed happen at all scales of `smoothness'. The applicability of the result is, naturally, depending on the growth of $\alpha_{\ell}$ in $\ell$, though, generically, one would not expect this to be badly behaved: a priori, there is no good reason to expect that the row of a matrix happens in any way to be the linear combination of singular vectors associated to large singular values -- though, naturally, this can happen (for example, if $A$ has one very large entry on the diagonal).

\section{Proofs}
\subsection{Proof of Theorem 1}
\begin{proof} To simplify exposition, we introduce the residual
$$ r_k = x_k - x.$$
Plugging in, we obtain that if the $i-$th equation is chosen, then
\begin{align*}
x + r_{k+1} &= x_{k+1} = x_k + \frac{ b_i - \left\langle a_i, r_k\right\rangle}{\|a_i\|^2} a_i \\
&=x +r_k + \frac{ b_i - \left\langle a_i, x + r_k\right\rangle}{\|a_i\|^2} a_i \\
&= x + r_k  - \frac{  \left\langle a_i,  r_k\right\rangle}{\|a_i\|^2} a_i + \left( \frac{ b_i - \left\langle a_i, x\right\rangle}{\|a_i\|^2} a_i\right).
\end{align*}
Since $x$ is the exact solution, we have $b_i - \left\langle a_i, x\right\rangle=0$ and
$$ r_{k+1} = r_k - \frac{  \left\langle a_i,  r_k\right\rangle}{\|a_i\|^2} a_i.$$
Recalling that the $i-$th row is chosen with probability proportional to $\|a_i\|^2$,
\begin{align*}
 \mathbb{E}~ \|A r_{k+1} \|^2 &= \mathbb{E}_i~\left\|A\left(  r_k - \frac{  \left\langle a_i,  r_k\right\rangle}{\|a_i\|^2} a_i \right) \right\|^2 \\
 &= \sum_{i=1}^{n} \frac{\|a_i\|^2}{\|A\|_F^2}\left\|A  r_k - \frac{  \left\langle a_i,  r_k\right\rangle}{\|a_i\|^2} Aa_i  \right\|^2.
 \end{align*}
 This norm can be explicitly squared out as 
 \begin{align*}
  \left\|A  r_k - \frac{  \left\langle a_i,  r_k\right\rangle}{\|a_i\|^2} Aa_i  \right\|^2 =  \|A {r_k}\|^2 - 2 \frac{ \left\langle a_i,  r_k\right\rangle}{\|a_i\|^2} \left\langle A r_k, A a_i \right\rangle + \frac{ \left\langle a_i,  r_k\right\rangle^2}{\|a_i\|^4} \|Aa_i\|^2.
 \end{align*}
 This allows us to rewrite the summation as
 \begin{align*}
 \mathbb{E}~ \|A r_{k+1} \|^2 &=  \sum_{i=1}^{n} \frac{\|a_i\|^2}{\|A\|_F^2} \left( \|A {r_k}\|^2 - 2 \frac{ \left\langle a_i,  r_k\right\rangle}{\|a_i\|^2} \left\langle A r_k, A a_i \right\rangle + \frac{ \left\langle a_i,  r_k\right\rangle^2}{\|a_i\|^4} \|Aa_i\|^2 \right)\\
 &= \|A r_k\|^2 - \frac{2}{\|A\|_F^2} \sum_{i=1}^{n}  \left\langle a_i,  r_k\right\rangle  \left\langle A r_k , A a_i \right\rangle  + \frac{1}{\|A\|_F^2} \sum_{i=1}^{n}  \left\langle a_i,  r_k\right\rangle^2 \frac{\| A a_i\|^2}{\|a_i\|^2}.
  \end{align*}
We have
\begin{align*}
  \frac{2}{\|A\|_F^2} \sum_{i=1}^{n}  \left\langle a_i,  r_k\right\rangle  \left\langle A r_k , A a_i \right\rangle &=  \frac{2}{\|A\|_F^2} \sum_{i=1}^{n}  \left\langle a_i,  r_k\right\rangle  \left\langle A^T A r_k , a_i \right\rangle \\
  &= \frac{2}{\|A\|_F^2}  \left\langle A r_k, A A^T A r_k \right\rangle\\
  &= \frac{2}{\|A\|_F^2}  \| A^T A r_k \|^2.
\end{align*}
The last sum we bound from above via
$$  \frac{1}{\|A\|_F^2} \sum_{i=1}^{n}  \left\langle a_i,  r_k\right\rangle^2 \frac{\| A a_i\|^2}{\|a_i\|^2} \leq \frac{\max_{i} \|A a_i\|^2/\|a_i\|^2}{\|A\|_F^2} \|A r_k\|^2.$$
This results in the desired estimate. 

\end{proof}

\subsection{Proof of Theorem 2}

\begin{proof} We again reduce the problem to that of the study of the residual
$$ r_{k+1} = r_k - \frac{  \left\langle a_i,  r_k\right\rangle}{\|a_i\|^2} a_i.$$
When looking at integer powers, we observe that, by the same reasoning,
\begin{align*}
 \mathbb{E}~ \|A^{\ell} r_{k+1} \|^2 &= \mathbb{E}~\left\|A^{\ell}\left(  r_k - \frac{  \left\langle a_i,  r_k\right\rangle}{\|a_i\|^2} a_i \right) \right\|^2 \\
 &= \sum_{i=1}^{n} \frac{\|a_i\|^2}{\|A\|_F^2}\left\|A^{\ell}  r_k - \frac{  \left\langle a_i,  r_k\right\rangle}{\|a_i\|^2} A^{\ell} a_i  \right\|^2\\
 &=  \sum_{i=1}^{n} \frac{\|a_i\|^2}{\|A\|_F^2} \left( \|A^{\ell} {r_k}\|^2 - 2 \frac{ \left\langle a_i,  r_k\right\rangle}{\|a_i\|^2} \left\langle A^{\ell} r_k, A^{\ell} a_i \right\rangle + \frac{ \left\langle a_i,  r_k\right\rangle^2}{\|a_i\|^4} \|A^{\ell} a_i\|^2 \right)\\
 &= \|A^{\ell} r_k\|^2 - \frac{2}{\|A\|_F^2} \sum_{i=1}^{n}  \left\langle a_i,  r_k\right\rangle  \left\langle A^{\ell} r_k , A^{\ell} a_i \right\rangle  \\
 &+ \frac{1}{\|A\|_F^2} \sum_{i=1}^{n}  \left\langle a_i,  r_k\right\rangle^2 \frac{\| A^{\ell} a_i\|^2}{\|a_i\|^2}.
  \end{align*}
  The first term is easy and the third term can, as before, be bounded by
  $$ \frac{1}{\|A\|_F^2} \sum_{i=1}^{n}  \left\langle a_i,  r_k\right\rangle^2 \frac{\| A^{\ell} a_i\|^2}{\|a_i\|^2} \leq \frac{\alpha_{\ell}}{\|A\|_F^2} \sum_{i=1}^{n}  \left\langle a_i,  r_k\right\rangle^2 =  \frac{\alpha_{\ell}}{\|A\|_F^2}  \|A r_k\|^2.$$
  It remains to understand the second term: here, we can use the symmetry of the matrix to write
  \begin{align*}
   \frac{2}{\|A\|_F^2} \sum_{i=1}^{n}  \left\langle a_i,  r_k\right\rangle  \left\langle A^{\ell} r_k , A^{\ell} a_i \right\rangle &= 
  \frac{2}{\|A\|_F^2} \sum_{i=1}^{n}  \left\langle a_i,  r_k\right\rangle  \left\langle A^{2\ell} r_k , a_i \right\rangle \\
  &=   \frac{2}{\|A\|_F^2}  \left\langle A r_k, A^{2\ell+1} r_k  \right\rangle \\
  &=\frac{2}{\|A\|_F^2}  \left\langle A^{\ell+1} r_k, A^{\ell+1} r_k  \right\rangle \\
  &= \frac{2}{\|A\|_F^2} \| A^{\ell+1} r_k\|^2.
  \end{align*}
\end{proof}

\end{document}